\documentclass{article}
\usepackage{latexsym}
\usepackage{graphics}
\usepackage{amsfonts}
\usepackage{graphicx}
\usepackage{epstopdf}
\newtheorem{thm}{Theorem}

\newtheorem{lem}{Lemma}
\newtheorem{defn}{Definition}

\begin{document}

\title{Tree Diagrams for String Links}
\author{Blake Mellor\\
			Mathematics Department\\
			Loyola Marymount University\\
			Los Angeles, CA  90045-2659\\
   {\it  bmellor@lmu.edu}}
\date{}
\maketitle

\vspace{.5in}

\begin{abstract}
In previous work, the author defined the intersection graph of a chord diagram
associated with a string link (as in the theory of finite type invariants).  In this paper,
we classify the trees which can be obtained as intersection graphs of string link diagrams.
\end{abstract}

\section{Introduction} \label{S:intro}

Intersection graphs of various kinds have an important role in graph theory \cite{mm}.  Recently,
intersection graph have also proved useful in topology, in the study of finite type knot invariants.  The
theory of finite type invariants allows us to interpret many important knot and link invariants in
purely combinatorial terms, as functionals on spaces of chord diagrams.  For knots, there is an obvious
intersection graph associated with these diagrams, first studied (in this context) by Chmutov, Duzhin
and Lando \cite{cdl}.  In many cases, these graphs contain all the relevant information for the
functionals coming from knot invariants \cite{cdl, me}.  Some of the most important knot invariants,
such as the Conway, Jones, Homfly and Kauffman polynomials, can be interpreted in terms of these
intersection graphs \cite{bg, me2}.

The machinery of finite type invariants and chord diagrams is easily extended to other topological
objects, such as links, braids and string links.  Recently, the author extended the notion of the
intersection graph of the chord diagram to string links as well \cite{me3}, presenting the
possibility of reinterpreting invariants of string links in the same way we have reinterpreted knot
invariants.  An important part of this program is to understand which graphs can appear as intersection
graphs in this context.  For knots, this question was answered by (among others) Bouchet \cite{bo}; for
string links, however, the class of graphs is much larger.  As special cases, it contains the set of
intersection graphs for knot diagrams (known to graph theorists as {\it circle graphs} \cite{mm}) and the
set of {\it permutation graphs}, both of which have been studied extensively.  This paper begins the
process of classifying these intersection graphs with a very modest goal: to describe the (directed)
trees which can arise as intersection graphs of chord diagrams for string links.

We will not discuss the background of finite type invariants, instead looking at the
relationship between the intersection graph and the chord diagram from a purely
combinatorial viewpoint.  For a discussion of how these diagrams arise in the theory of
finite type invariants, see \cite{bn}.  In Section~\ref{S:prelim} we will review the
definitions of chord diagrams and intersection graphs for string links.  In
Section~\ref{S:classify} we will classify the trees which arise as intersection graphs.

{\it Acknowledgement:}  The author thanks Loyola Marymount University
for supporting this work via a Summer Research Grant in 2003.

\section{Preliminaries} \label{S:prelim}

\subsection{Chord Diagrams} \label{SS:chord}

We begin by defining what we mean by a {\it chord diagram}.  Since we are only considering
chord diagrams which arise from string links, we will simply refer to these as
chord diagrams, but the reader should be aware that these differ somewhat from the more
usual chord diagrams for knots.

\begin{defn}
A {\bf chord diagram of degree n with k components} is a disjoint union of k oriented
line segments (called the {\bf components} of the diagram), together with $n$ chords (unoriented line
segments with endpoints on the components), such that all of the $2n$ endpoints of the chords are
distinct.  The diagram is determined by the orders of the endpoints on each component.
\end{defn}

In the theory of finite type knot and link invariants, these diagrams are usually organized into a
graded vector space, considered modulo certain relations.  For the purposes of this paper, however, this
additional structure is irrelevant.

\subsection{Intersection Graphs} \label{SS:IG}

The essential value of the intersection graph for knots (in which the chord diagram
consists of chords in a bounding circle) is that it can detect when the order of two
endpoints for different chords along the bounding circle is switched, since this changes the
pair of chords from intersecting to non-intersecting or vice-versa.  For chord
diagrams for string links, the existence of a "bottom" and "top" for each component allows us
to give a linear (rather than cyclic) ordering to the endpoints of the chords on each
component, and so the notion of one endpoint being "below" another is well-defined.  We want
our intersection graphs to detect when this order is reversed.

\begin{defn} \cite{me3} \label{D:IGstring}
Let $D$ be a chord diagram with $k$ components (oriented line
segments, colored from 1 to $k$) and $n$ chords.  The {\it intersection graph}
$\Gamma(D)$ is the labeled, directed multigraph such that:
\begin{itemize}
    \item $\Gamma(D)$ has a vertex for each chord of $D$.  Each vertex is labeled
by an unordered pair $\{i,j\}$, where $i$ and $j$ are the labels of the components
on which the endpoints of the chord lie.
    \item There is a directed edge from a vertex $v_1$ to a vertex $v_2$ for each
pair $(e_1, e_2)$ where $e_1$ is an endpoint of the chord associated to $v_1$,
$e_2$ is an endpoint of the chord associated to $v_2$, $e_1$ and $e_2$ lie on the
same component of $D$, and the orientation of the component runs from $e_1$ to
$e_2$ (so if the components are all oriented upwards, $e_1$ is below $e_2$).  We
count these edges "mod 2", so if two vertices are connected by two
directed edges with the same direction, the edges cancel each other.  If two vertices
are connected by a directed edge in each direction, we will simply connect them
by an undirected edge.
\end{itemize}
\end{defn}

\begin{figure} [h]
    $$\includegraphics{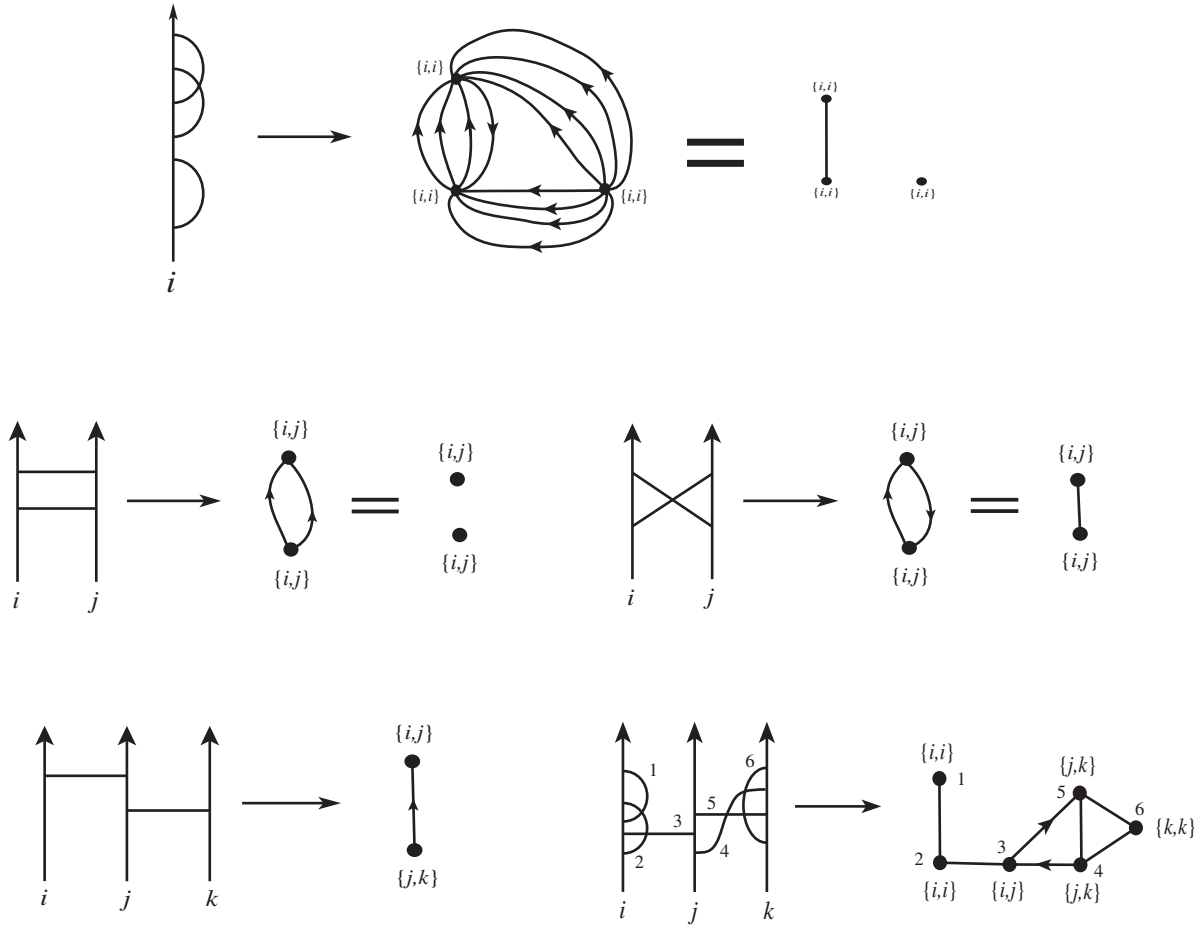}$$
    \caption{Examples of intersection graphs for string links.  In the last example the chords and vertices are numbered to clarify the bijection between them.} \label{F:IGstring}
    \end{figure}

Examples of chord diagrams and their associated intersection graphs are given in
Figure~\ref{F:IGstring}.
Note that when the two chords have both endpoints on the same component $i$, our
definition of intersection graph corresponds to the usual intersection graph for
knots.  Our definition also matches our intuition in the case of chord diagrams
of two components, as shown in Figure~\ref{F:IGstring}.

Note also that the total number of directed edges between a vertex $v$ labeled
$\{i,j\}$ and a vertex $w$ labeled $\{l,m\}$ is given by the sum of the number
of occurrences of $i$ in $\{l,m\}$ and the number of occurrences of $j$ in
$\{l,m\}$.  In particular, if a vertex $v$ has a label $\{i,i\}$, this number
will be even (0, 2 or 4).  Since we count directed edges modulo 2, this implies
there is an (uncancelled) directed edge from $v$ to another vertex $w$ if and only if
there is also an (uncancelled) directed edge from $w$ to $v$.  We will say that labeled
directed multigraphs which have this property are {\it semisymmetric}.

\begin{defn} \label{D:semisym}
A directed multigraph G, with each vertex labeled by a pair \{i,j\}, is {\bf
semisymmetric} if for every vertex v labeled \{i,i\}, and any other vertex w, there is a
directed edge from v to w if and only if there is a directed edge from w to v.
\end{defn}

\section{Intersection Graphs and Trees} \label{S:classify}

In this section we will give a characterization of the trees which arise as intersection
graphs.  For our purposes, a directed graph is a {\it tree} if it is an {\it
undirected} tree - i.e. the graph is connected, and it has no cycles even if the orientations
of the edges are ignored.  We will first consider chord diagrams with two components, and
then the (easier) case of diagrams with more than two components.

\subsection{2-component diagrams} \label{SS:2comp}

In intersection graphs for diagrams with two components, labeled 1 and 2, there are three possible labels
on the vertices:  \{1,1\}, \{2,2\} and \{1,2\}.  To begin with, we will consider {\it marked}, rather
than {\it labeled} intersection graphs.  A vertex will be called {\it marked} if it is labeled \{1,2\},
and {\it unmarked} if it is labeled \{1,1\} or \{2,2\}.  In other words, a vertex is marked if the
corresponding chord has its endpoints on different components, and unmarked if the chord has both
endpoints on the same component.

\begin{figure} [h]
    $$\includegraphics{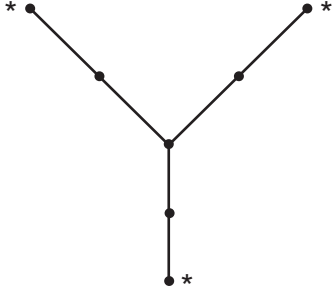}$$
    \caption{A marked graph which is not an intersection graph} \label{F:nonIG}
    \end{figure}

The first important observation is that not all marked trees are intersection graphs
(otherwise the characterization would be rather trivial).  Figure~\ref{F:nonIG} gives an
example of a marked tree which is not an intersection graph for any chord diagram on two
components.  The marked vertices are labeled by an asterisk (*).

To determine which marked trees {\it are} intersection graphs, we need to introduce the
notion of {\it light} and {\it heavy boughs} of a tree.

\begin{defn} \label{D:bough}
A {\bf bough} of a vertex v in a marked tree T is a connected component of $T\backslash v$
(the graph which results from removing the vertex v and all edges incident to v).  The bough
is called {\bf light} if it contains at most one marked vertex and this vertex (if present)
is adjacent to v in T.  Otherwise the bough is {\bf heavy}.
\end{defn}

The main theorem of this section characterizes tree intersection graphs via boughs.

\begin{thm} \label{T:IGtree}
A marked tree T (i.e. T may or may not have marked vertices) is an intersection graph for a chord
diagram on two components if and only if any vertex in T has at most 2 heavy boughs.
\end{thm}

Before we prove this theorem, we require two lemmas.  Our first lemma relates the notion
of a light bough to the notion of a {\it share} of a chord diagram, introduced in \cite{cdl}.

\begin{defn} \label{D:share}
A {\bf share} of a chord diagram D is a subset S of the set of chords in D and two arcs A
and B on the boundary components of D (A and B may be on the same or different components)
such that every chord in S has both endpoints in $A \cup B$ and {\bf no} other chord in D
has an endpoint in $A \cup B$.
\end{defn}

\begin{lem} \label{L:share}
Assume T is a marked tree which is the intersection graph for a chord diagram D, and v is a
vertex in T.  Further assume that $T$ has at least one marked vertex.  Then a bough of v is
light if and only if the corresponding chords are a share in D.
\end{lem}
{\sc Proof:}  Recall that trees are {\it connected} graphs.  We will first assume
that the chords corresponding to the vertices in a bough $H$ of $v$ form a share $S$
in $D$, and show that $H$ is light.  We will show that any marked vertex in $H$ must
be adjacent to $v$.  Working by contradiction, assume $w$ is a marked vertex in $H$ which is
not adjacent to $v$.  Since $S$ is a share, all endpoints of the chords in $S$ lie on two
arcs $A$ and $B$.  Since the chord corresponding to $w$ (which we will also call $w$) has an
endpoint on each component of the chord diagram, we may assume that $A$ and $B$ lie on
different components, and $w$ has one endpoint in each.  Since $v$ is {\it not} in $S$, it
has no endpoints in $A$ or $B$.  Since $v$ and $w$ are not adjacent, the endpoints of $v$
must be either both above or both below the endpoints of $w$.  In this case, the endpoints
of $v$ are both above or both below {\it any} endpoint in $A \cup B$, so $v$ is not
adjacent to any vertex in $H$ (since the directed edges in an intersection graph are counted "mod 2"). 
But since $T$ is connected, and $H$ is a component of $T\backslash v$, this is a contradiction, and we
conclude that any marked vertex in $H$ is adjacent to $v$.

We still need to show that $H$ has at most one marked vertex.  Assume $H$ has two marked
vertices, $w$ and $u$.  By the argument above, these are both adjacent to $v$, so there is
a path $wvu$ in $T$.  But, since $H$ is a connected component of $T\backslash v$, there is also a
path from $w$ to $u$ in $H$ which does not contain $v$.  Joining these two paths gives us a
cycle in $T$, which contradicts the assumption that $T$ is a tree, so $H$ can have at most
one marked vertex.  This proves $H$ is a light bough of $v$.

Now we assume that $H$ is a light bough of $v$, and $S$ is the set of chords corresponding to the
vertices of $H$.  Since $H$ is connected, and $T$ is a tree, there is a unique vertex $w$ in
$H$ which is adjacent to $v$ (otherwise, we will have a cycle as in the last paragraph). 
$w$ may or may not be marked - in either case, no other vertices of $H$ are marked.  Our
goal will be to show that we can find minimal arcs $A$ and $B$ containing all the endpoints
of chords in $S$ (but not the endpoints of $v$), and these will be the arcs required for $S$
to be a share.  The first step is to show that there exists some pair of arcs containing all
the endpoints of chords in $S$, but not the endpoints of $v$.

\begin{figure} [h]
    $$\includegraphics{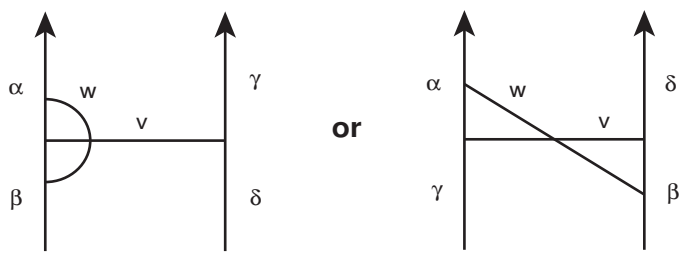}$$
    \caption{$v$ is marked} \label{F:lemma1}
    \end{figure}

We first consider the case when $v$ is marked.  Then $v$ divides each of the two boundary
components of the diagram into two pieces (giving a total of four segments), as shown in
Figure~\ref{F:lemma1}.
$w$ connects two of these four pieces, $\alpha$ and $\beta$.  Since no other chords in $S$
intersect $v$, and all other chords in $S$ are unmarked, the other chords in $S$ must have
both endpoints in $\alpha$ or both endpoints in $\beta$.  So $\alpha$ and $\beta$ are arcs
containing all the endpoints of $S$, but not the endpoints of $v$.

\begin{figure} [h]
    $$\includegraphics{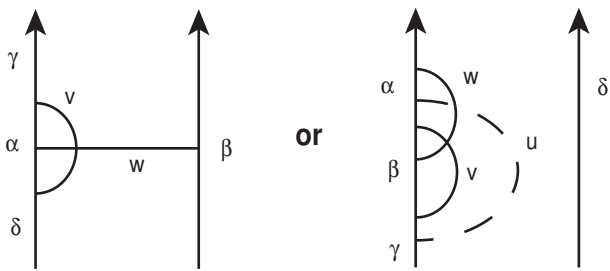}$$
    \caption{$v$ is unmarked} \label{F:lemma2}
    \end{figure}

If $v$ is unmarked, it again divides the two boundary components into four pieces, three on
one component and one on the other.  Once again, $w$ connects two of these pieces, $\alpha$
and $\beta$ as shown in Figure~\ref{F:lemma2}.
If $w$ is marked, then once again all chords of $S$ must have both endpoints in $\alpha$ or
both in $\beta$, and we can take these as the desired arcs.  If $w$ is also unmarked, there
may be a chord $u$ in $S$ which connects $\alpha$ to one of the other arcs, as shown in
Figure \ref{F:lemma2}.  But this means there are {\it no} marked chords in the diagram:  any
marked chord would have to be connected to $v$ (since the tree $T$ is connected) by a path of
chords which would at some point cross $u$, but this would mean $S$ contains a marked chord,
which is impossible.  Since our lemma assumes that $T$ has at least one marked vertex, we
can discard this case.  If there is no such chord $u$, then we can take the arcs $\alpha$
and $\beta$ as before.

Now that we know there exists at least one pair of disjoint arcs on the boundary components
of the chord diagram which contain all the endpoints of $S$ but not the endpoints of $v$, let
$A$ and $B$ be a {\it minimal} pair of such arcs - minimal in the sense of containing as few
endpoints of chords in $T$ as possible.  Now say we have a chord $u$ which is not in
$S$.  First assume that $u$ has one endpoint in $A \cup B$; without loss of generality,
we may assume it lies in $A$.  Since $A \cup B$ is minimal, there are endpoints of chords in
$S$ both above and below $u$ in $A$ (otherwise we could shrink $A$ to exclude the endpoint
of $u$, giving a smaller arc).  Since $H$ is connected, there is a path of chords connecting the
endpoints above $u$ to those below $u$.  If the other endpoint of $u$ is not in $A$, it must
cross some chord in this path, but this means that $u$ is in $S$, a contradiction.  So the other
endpoint of $u$ must also be in $A$.  By the same argument, any chord intersecting $u$ must have
at least one, and hence both, endpoints in $A$; in fact, all chords corresponding to vertices in
the connected component of $T$ containing $u$ must have both endpoints in $A$.  But since $T$ is
connected, this means that, in particular, $v$ has both endpoints in $A$, a contradiction.  We
conclude that {\it only} chords in $S$ have endpoints in $A \cup B$, which shows that $S$ is a
share, and completes the proof of the lemma.  $\Box$

\begin{figure} [h]
    $$\includegraphics{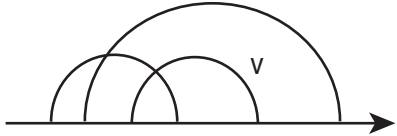}$$
    \caption{A light bough which is not a share.} \label{F:notshare}
    \end{figure}

Note that if the tree $T$ does {\it not} have any marked vertices, a vertex $v$ may have a
light bough which does not correspond to a share, as in Figure~\ref{F:notshare}.
However, in this case we are only looking at diagrams on a single component.  Bar-Natan
\cite{bn} showed that this case is equivalent (modulo the 4-term relation) to chord diagrams
on knots, and if we change the bounding components to a circle the boughs {\it are} shares
\cite{cdl}.  Since Chmutov {\it et al.} \cite{cdl} have already shown that all unmarked
trees are intersection graphs, in this section we restrict our attention to trees with at
least one marked vertex.

We now prove a lemma about the structure of trees where every vertex has at most two heavy
boughs.  We introduce the notion of the {\it spine} of a tree:

\begin{defn} \label{D:spine}
Say that $T$ is a marked tree.  A {\bf spine} of $T$ is a path in $T$ of maximal length such that
the first and last vertices in the path are marked.  If $T$ has a single marked vertex, that
vertex is the spine; if it has no marked vertices the spine is empty.
\end{defn}

\begin{figure} [h]
    $$\includegraphics{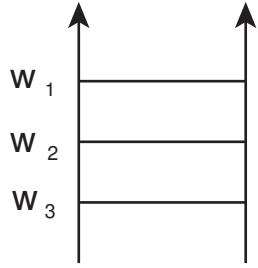}$$
    \caption{Parallel marked chords.} \label{F:parallel}
    \end{figure}

\begin{lem} \label{L:spine}
If $T$ is a marked tree where every vertex has at most two heavy boughs, and $S$ is a spine for
$T$, then every marked vertex in $T$ which is not in $S$ is adjacent to a vertex in the interior
(i.e. not one of the endpoints) of $S$.
\end{lem}
{\sc Proof:}  Let $v$ be a marked vertex in $T$ which is not in $S$.  Since $T$ is connected,
there is a path $P$ connecting $v$ to some vertex $w$ in $S$ (such that no other vertex on $P$ is
in $S$).  Since $S$ has maximal length, $w$ must be an interior vertex of $S$, or the path $P
\cup S$ would be longer than $S$.  Also, $P$ must be no longer than the path from $w$ to either
endpoint of $S$, or we could combine $P$ with one of these paths to get a path between marked
vertices longer than $S$.  So if $v$ is not adjacent to $w$, neither are either of the endpoints
of $S$.  But this means that $w$ has three heavy boughs.  So $v$ must be adjacent to $w$, and we
can conclude that every marked vertex of $T$ which is not in $S$ is adjacent to an interior
vertex of $S$. $\Box$
\\
\\
\noindent{\sc Proof of Theorem \ref{T:IGtree}:}  We first assume that $T$ is an intersection
graph, and show that any vertex $v$ has at most two heavy boughs.  We work by contradiction -
assume $v$ is a vertex with (at least) three heavy boughs; call these boughs $C_1,\ C_2,\ C_3$. 
So each $C_i$ contains a marked vertex $w_i$ which is {\it not} adjacent to $v$.  Since
$w_1,\ w_2,\ w_3$ are non-adjacent marked vertices, the corresponding chords must be parallel, as
shown in Figure~\ref{F:parallel}.
Without loss of generality, we order the chords $w_1,\ w_2,\ w_3$ from top to bottom.  Since
$T$ is a tree, there is a unique path from $w_1$ to $w_3$, passing through $v$, and no
chord in this path is adjacent to $w_2$ (since the only chord in the path adjacent to any
chord in $C_2$ is $v$ itself, and $v$ is not adjacent to $w_2$).  But since $w_2$ separates
$w_1$ and $w_3$, {\it some} chord in the path must cross $w_2$, which is a contradiction. 
So any chord $v$ can have at most two heavy boughs.

\begin{figure} [h]
    $$\includegraphics{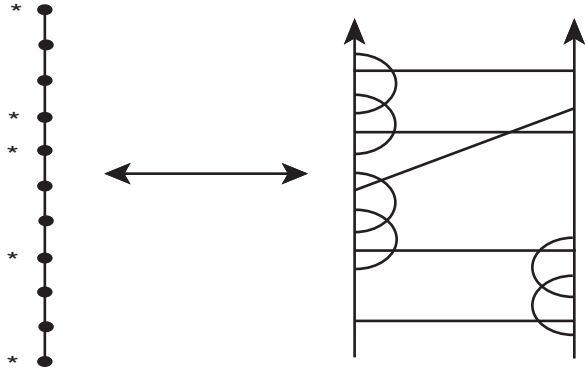}$$
    \caption{Chord diagram for a path.} \label{F:spinediagram}
    \end{figure}

\begin{figure} [h]
    $$\includegraphics{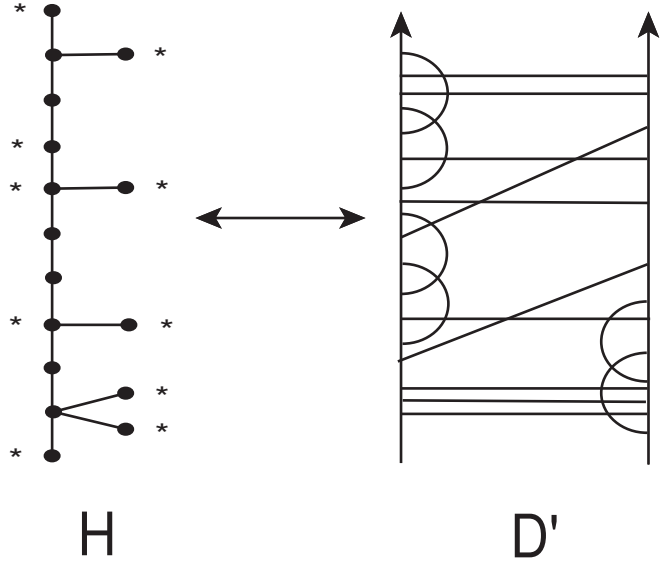}$$
    \caption{Expanding $D$ to $D'$} \label{F:ribs}
    \end{figure}

\begin{figure} [h]
    $$\includegraphics{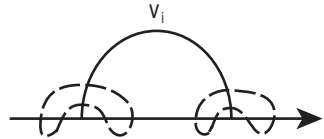}$$
    \caption{Chord diagram for $C_i$.} \label{F:barbell}
    \end{figure}

\begin{figure} [h]
    $$\includegraphics{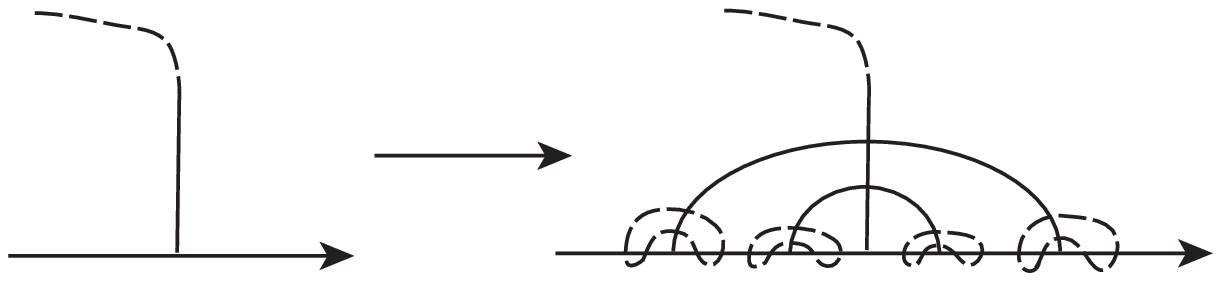}$$
    \caption{Adding "barbells" to a chord diagram} \label{F:addingbarbells}
    \end{figure}

Now we will assume that any vertex $v$ has at most two heavy boughs, and show that $T$ is an
intersection graph.  Let $S$ be a spine for $T$.  $S$ is just a path of marked and unmarked
vertices, so $S$ is the intersection graph for a chord diagram $D$, as shown in
Figure~\ref{F:spinediagram}.
By Lemma \ref{L:spine}, any marked vertex of $T$ is adjacent to an interior vertex of $S$; let
$H$ be the union of $S$ and the other marked vertices of $T$.  Then we can add chords to $D$ to
produce a new diagram $D'$ with intersection graph $H$, as shown in Figure~\ref{F:ribs}.
All remaining vertices of $T$ are unmarked.  These vertices can be divided into connected
components $C_i$, each of which is adjacent to one vertex of $H$.  Let $v_i$ denote the vertex
of $C_i$ which is adjacent to $H$.  By \cite{cdl}, $C_i$ is the intersection graph for a
"barbell" chord diagram consisting of the chord corresponding to $v_i$ and a share at either end,
as shown in Figure~\ref{F:barbell}.
These diagrams can be attached to the desired chords of $D'$ (and any number can be attached to a
single chord of $D'$) to obtain a chord diagram with intersection graph $T$, as in
Figure~\ref{F:addingbarbells}.
So $T$ is an intersection graph, completing the proof.  $\Box$

\begin{figure} [h]
    $$\includegraphics{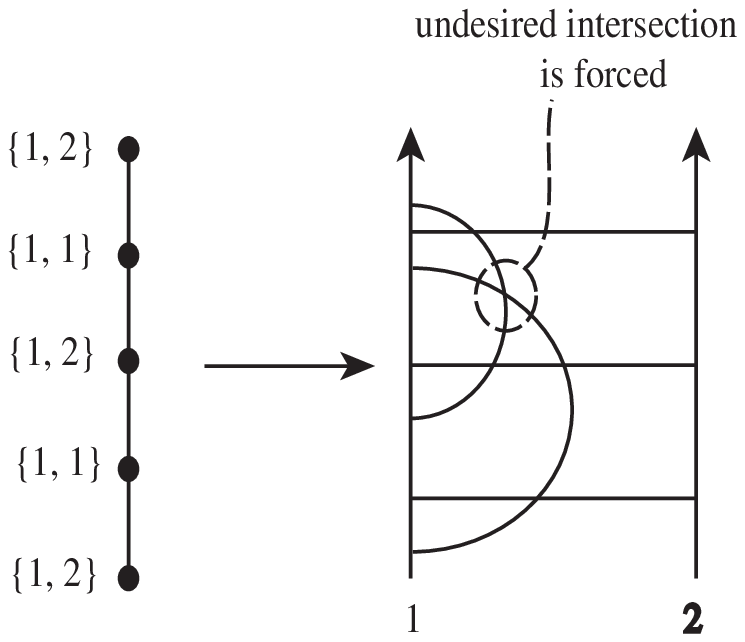}$$
    \caption{A labeled tree which is not an intersection graph.} \label{F:badlabel}
    \end{figure}

We now turn our attention from {\it marked} trees to {\it labeled} trees, in which the marked
vertices are labeled \{1,2\}, and the unmarked vertices are labeled either \{1,1\} or \{2,2\},
and ask when a labeled tree is an intersection graph.  There are two obvious requirements:  the
underlying marked tree must be an intersection graph (so the tree must satisfy the condition of
Theorem \ref{T:IGtree}), and adjacent unmarked vertices must have the same label (since chords on
different components of a chord diagram cannot intersect).  However, these conditions are not
sufficient, as shown by the example in Figure~\ref{F:badlabel}.
We need to add one more condition, which is given by the following theorem:

\begin{thm} \label{T:IGlabel}
A labeled tree $T$ is an intersection graph for a chord diagram on two components if and only if it
satisfies the following three conditions:
\begin{itemize}
     \item  The underlying marked tree is an intersection graph.
	    \item  Adjacent unmarked vertices have the same label.
     \item  If $S$ is a spine of $T$, and $v$ and $w$ are unmarked vertices in $S$, then the
labels on $v$ and $w$ are the same if and only if there is an even number of marked vertices
between them along $S$.
\end{itemize}
\end{thm}
{\sc Proof:}  We first prove the sufficiency.  If a labeled tree $T$ satisfies the three
conditions, then we can construct a chord diagram with intersection graph $T$ as in Theorem
\ref{T:IGtree}.  In particular, the third condition allows us to construct a chord diagram for
the spine as in Figure \ref{F:spinediagram}.

As for the necessity, it is clear that the first two conditions are necessary for $T$ to be an
intersection graph.  Assume $T$ is an intersection graph for a chord diagram $D$.  Let $S$ be a
spine of $T$, and $v$ and $w$ be two unmarked vertices along $S$.  If there are {\it no} marked
vertices between $v$ and $w$, then they are connected by a path of unmarked vertices corresponding
to chords in the chord diagram with both endpoints on the same component.  So all of these chords
must be on the same component, and the vertices must have the same label.  Now assume there are
marked vertices between $v$ and $w$.

\begin{figure} [h]
    $$\includegraphics{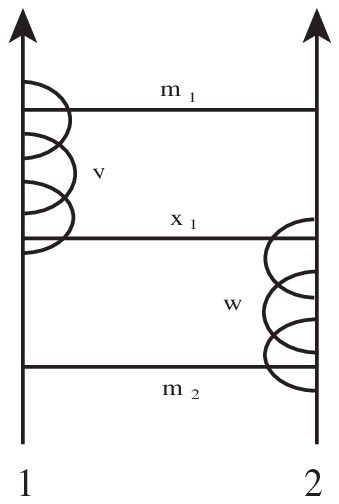}$$
    \caption{$v$ and $w$ separated by one marked chord.} \label{F:labelbase}
    \end{figure}

We will first consider the case when all of the marked vertices between $v$ and $w$ along $S$ are
consecutive along $S$.  Since the endpoints of the path $S$ are marked vertices, there is a path
$m_1U_1vU_2MU_3wU_4m_2$ in $S$, where $m_1$ and $m_2$ are marked vertices, each $U_i$ is a path
of unmarked vertices, and $M = x_1x_2...x_k$ is a path of $k$ marked vertices.  Since $T$ is a
tree, $m_1$ and $m_2$ are not adjacent to each other or to any of the vertices in $M$.  If $k =
1$, we have the picture in Figure~\ref{F:labelbase} (assuming, without loss of generality, that
$U_1vU_2$ is on component 1).
Since the first chord in $U_3$ does not intersect any chord in $U_1vU_2$, it cannot have an
endpoint on component 1 between $m_1$ and $x_1$.  However, since it is intersects $x_1$, it
must have an endpoint on either side of $x_1$, and since it does not intersect $m_1$, one of
these endpoints must lie between $x_1$ and $m_1$.  So this chord must have an endpoint on
component 2.  Since it is unmarked, it will have both endpoints on component 2, so the entire
path $U_3wU_4$ is on component 2, and $v$ and $w$ have different labels.

\begin{figure} [h]
    $$\includegraphics{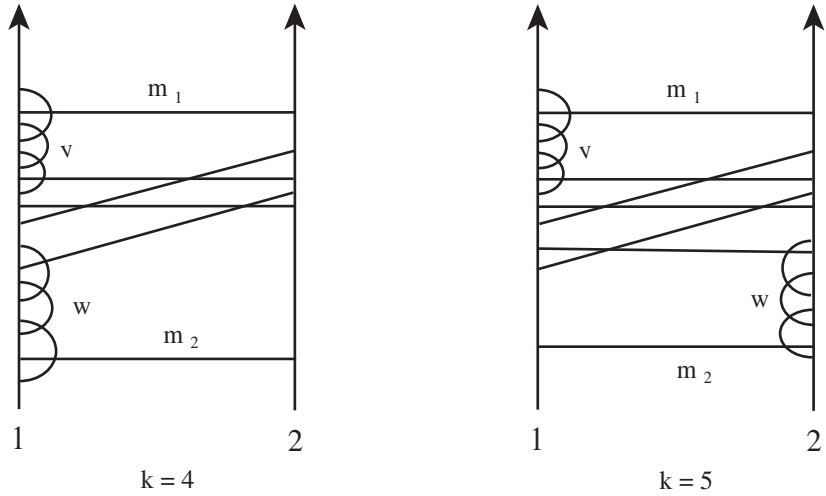}$$
    \caption{$v$ and $w$ separated by odd or even number of marked chords.} \label{F:labelgeneral}
    \end{figure}

If $k > 1$, the marked chords will cross each other as in Figure~\ref{F:labelgeneral}, which shows
examples when $k$ is even or odd.
This arrangement is forced by the fact that each of the marked vertices is only adjacent to the ones
immediately before and after it in the path $M$.  The first chord of $U_3$ must intersect $x_k$ on the
same side of $x_{k-1}$ as $m_2$.  If $k$ is odd, this means that the chords in $U_3wU_4$ will be on
component 2; if $k$ is even, they will lie on component 1.  So $v$ and $w$ will have the same label if
and only if $k$ is even.

If the marked chords between $v$ and $w$ are not all consecutive, we can apply the argument above to
each consecutive grouping.  The labels of $v$ and $w$ will be the same if and only if the labels switch
an even number of times, i.e. if there are an even number of stretches of marked vertices of odd length,
which will happen exactly when the total number of marked chords in between $v$ and $w$ is even.  This
completes the proof. $\Box$

\subsection{Diagrams with more than 2 components.} \label{SS:>2comp}

Surprisingly, the directed trees which are intersection graphs for string link chord
diagrams with more than 2 components is somewhat simpler than for diagrams with two
components - we no longer need to consider heavy boughs.  We will consider {\it connected} chord
diagrams, meaning that we cannot separate the components into two nonempty sets so that all chords are
between components in the same set.  Our goal in this section is to prove the following theorem:

\begin{thm} \label{T:>2comp}
Let T be a connected, labeled directed tree with n colors (so each vertex has a label \{i, j\}, where $1
\leq i, j\leq n$).  Let $m_{i,j}$ denote the number of vertices with label \{i, j\}.  Then T is an
intersection graph for a connected chord diagram on n components if and only if the following
conditions are met (possibly after relabeling the tree by a permutation of 1,..., n):
\begin{enumerate}
   \item The labels of adjacent vertices must have at least one color in common.
   \item T is semisymmetric (see Definition \ref{D:semisym}).
   \item If v has label \{i,j\} and w has label \{i,k\}, where i, j and k are all distinct, then there is
a directed edge between v and w.
   \item $m_{i,j} = 0$ if $|i-j| > 1$.
   \item $m_{i, i+1} = 1$ for $2\leq i\leq n-2$, $m_{1,2}\geq 1$ and $m_{n-1,n} \geq 1$.
   \item No two marked vertices are connected by a path of undirected edges.
\end{enumerate}
\end{thm}

\noindent{\sc Remark:}  It is not intuitive that the theorem for many components should be simpler than that for two components - why could we not augment a 2-component chord diagram to a 3-component diagram by adding an additional component with a single chord connecting it to the previous diagram?  The problem is that if the original diagram had a heavy bough (and so two marked chords connected by a path of unmarked chords), then the new diagram would not be a tree diagram (since the new chord would be adjacent to both of the marked chords in the heavy bough, completing a cycle).  So heavy boughs are no longer permitted (this is one effect of Condition 6), and the diagrams are simpler.

\noindent{\sc Proof:}  We will maintain the convention of the previous section, and refer to a vertex as {\it marked} if it is
labeled by two separate colors and {\it unmarked} if it is labeled by only one color.  We will first show
that these conditions are all necessary.  Assume that a tree $T$ is an intersection graph for a connected
string link chord diagram $D$.  Condition 1 is obvious, since two chords cannot intersect (in the sense
of Definition \ref{D:IGstring}) unless they have at least one endpoint on the same component. 
Conditions 2 and 3 are consequences of the discussion at the end of Section \ref{SS:IG}.

To show the necessity of Condition 4, we will use the {\it connection graph} of the diagram $D$.

\begin{defn} \label{D:connection}
Given a string link chord diagram $D$ with $n$ components (labeled 1, 2,..., $n$), the {\bf connection
graph} of $D$ is the graph $C(D)$ with $n$ vertices (labeled 1, 2,..., $n$) where vertices $i$ and $j$
are adjacent in $C(D)$ if and only if there is a chord connecting components $i$ and $j$ in $D$.
\end{defn}

We first observe that if $T = \Gamma(D)$ is a tree, then so is $C(D)$.  If there is a cycle
$i_1i_2...i_k$ in $C(D)$, then there are chords in $D$ connecting components $i_1$ and $i_2$, $i_2$ and
$i_3$,..., $i_k$ and $i_1$.  So $T$ has vertices labeled $\{i_1,i_2\}$, $\{i_2,i_3\}$,...,
$\{i_{k-1},i_k\}$, $\{i_k,i_1\}$.  Call these vertices $v_1$,..., $v_k$.  Then, by Condition 3, each
$v_i$ is adjacent to $v_{i+1}$ (and $v_k$ is adjacent to $v_1$), so there is a cycle $v_1v_2...v_k$ in
$T$.  This is impossible, since $T$ is a tree, so $C(D)$ must also be a tree.

Since $D$ is a connected diagram, every vertex of $C(D)$ has degree at least one.  Moreover, no vertex
of $C(D)$ can have degree greater than two.  Otherwise, there would be a vertex $i$ in $C(D)$ adjacent
to three other vertices $j, k, l$.  So there would be three chords in $D$ with one endpoint on component
$i$ and the other on component $j, k, l$ (respectively).  This means there are corresponding vertices in
$T=\Gamma(D)$ with labels $\{i,j\}$, $\{i,k\}$ and $\{i,l\}$.  By Condition 3, these three vertices are
all adjacent to each other, meaning that $T$ will have a 3-cycle.  Since $T$ is a tree, this is
impossible.

So $C(D)$ is a connected tree where every vertex has degree at most two; this means that $C(D)$ is a
path.  Starting at one end of the path, we renumber the vertices (and hence the components of $D$) from
1 to $n$.  We will now have no chords between components $i$ and $j$ in $D$ if $|i-j| > 1$, which means
$m_{i,j} = 0$, giving us Condition 4.  Moreover, we have that $m_{i, i+1} \geq 1$ for $1\leq i\leq n-1$.

For the fifth condition, assume that $m_{i, i+1} > 1$ for some $i$ with $2\leq i\leq n-2$.  So $T$ will
have at least two vertices $v$ and $w$ labeled $\{i,i+1\}$.  Since $m_{i-1,i} \geq 1$ and $m_{i+1,i+2}
\geq 1$, we also have vertices $s$ and $t$ labeled $\{i-1,i\}$ and $\{i+1,i+2\}$, respectively.  By
Condition 3, $s$ and $t$ are both adjacent to both $v$ and $w$, so we have a cycle $svtw$ in $T$.  Since
$T$ is a tree, this is impossible, so $m_{i, i+1} = 1$ for $2\leq i\leq n-2$, giving us Condition 5.

For the final condition, assume that two marked vertices $v$ and $w$ are connected by an undirected
path $vu_1u_2...u_rw$ ($r \geq 1$ by Condition 3).  Then all the unmarked vertices in the path have the
same label, say $\{i,i\}$, which is one of the labels of $v$ and $w$.  If $v$ and $w$ share only one
label, there is also a directed edge between them, giving a cycle in the graph $T$.  Since $T$ is a
tree, this is impossible, so $v$ and $w$ must share both their labels, $\{i,i\pm 1\}$.  By Condition 5,
this label is either $\{1,2\}$ or $\{n-1,n\}$.  Without loss of generality, assume the label is
$\{1,2\}$.  Then there is a vertex $x$ labeled $\{2,3\}$ which is adjacent to both $v$ and $w$ by
directed edges, giving a cycle $vu_1u_2...u_rwxv$ in $T$.  Again, since $T$ is a tree, this is
impossible, which proves the final condition.  So any tree that is an intersection graph will satisfy the
six conditions of the theorem.

\begin{figure} [h]
    $$\includegraphics{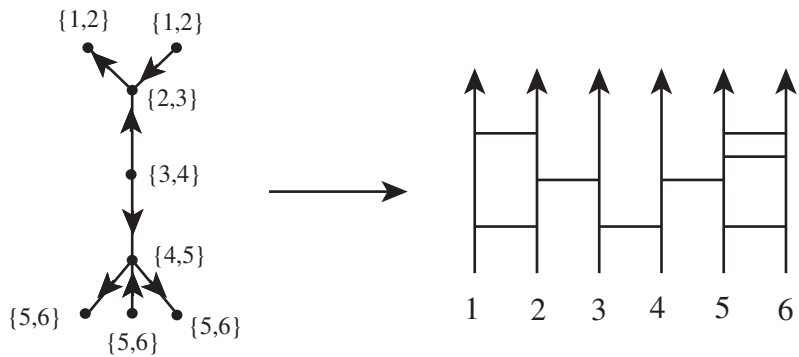}$$
    \caption{A chord diagram with intersection graph $S$} \label{F:makediagram}
    \end{figure}

It remains to show that these six conditions are sufficient.  Let $T$ be a connected, labeled directed
tree with $n$ colors which satisfies the conditions of the theorem.  If necessary, relabel the vertices
so Conditions 4 and 5 are satisfied.  We first consider the subgraph $S$ of $T$ induced by the marked
vertices; $S$ will still be connected by Conditions 3 and 5.  We will first construct a chord diagram
with intersection graph $S$.  Begin with an empty chord diagram with $n$ components (all oriented
upwards).  Add $m_{1,2}$ parallel chords between components 1 and 2.  (Note that no two vertices labeled
$\{1,2\}$ are adjacent, or they would form a 3-cycle with the unique, by Condition 5, vertex with lable
$\{2,3\}$.)  If $v$ is the unique vertex in $S$ with label $\{2,3\}$, there are $k$ vertices labeled
$\{1,2\}$ with an edge directed {\it towards} $v$, and $m_{1,2}-k$ such vertices with an edge directed
{\it away} from $v$.  Place a chord between components 2 and 3 so that the endpoint on component 2 is
above $k$ of the chords between components 1 and 2, and below the rest.  We now place a chord between
components 3 and 4, which is above or below the chord between components 2 and 3 depending on the
direction of the edge in $S$.  This corresponds to the unique (by Condition 5) vertex in $S$ with
label $\{3,4\}$.  Continue in this way through all the components, finally putting $m_{n-1,n}$ parallel
chords between components $n-1$ and $n$.  The result is a chord diagram with intersection graph $S$. 
See Figure~\ref{F:makediagram} for an example.  
Now we simply need to add chords corresponding to the unmarked vertices in $T$.  By Condition 6, no
marked vertices are connected by paths of unmarked vertices.  So the unmarked "branches" of the tree
correspond to "barbell" chord diagrams, as in Figure \ref{F:barbell}, and these diagrams can be attached
to the chord diagram for $S$ as in Figure \ref{F:addingbarbells}.  The result is a chord diagram with
intersection graph $T$, which completes the proof of Theorem \ref{T:>2comp}. $\Box$

\end{document}